\documentclass[12pt]{article}
\usepackage{amssymb}

\usepackage{amsmath,latexsym,amsfonts}
\usepackage{graphicx}
\usepackage{srcltx}
\usepackage{amscd}

\newtheorem{lemma}{Lemma}[section]
\newtheorem{coro}[lemma]{Corollary}
\newtheorem{prop}[lemma]{Proposition}
\newtheorem{thm}[lemma]{Theorem}
\newtheorem{defn}[lemma]{Definition}
\newtheorem{rem}[lemma]{Remark}
\makeatletter\@addtoreset{equation}{section}
\renewcommand\theequation{\thesection.\@arabic\c@equation}

\begin{document}
\begin{center}
{\LARGE   $L^2$-transverse conformal Killing forms on complete foliated manifolds}

 \renewcommand{\thefootnote}{}
\footnote{2000 \textit {Mathematics Subject Classification.}
53C12, 53C27, 57R30}\footnote{\textit{Key words and phrases.}
Transverse Killing form, Transverse conformal Killing form}
\renewcommand{\thefootnote}{\arabic{footnote}}
\setcounter{footnote}{0}

\vspace{1 cm} {\large Seoung Dal Jung and Huili Liu}
\end{center}
\vspace{0.5cm}

{\bf Abstract.} In this article, we study the $L^2$-transverse conformal Killing forms on complete foliated Riemannian manifolds and prove some vanishing theorems. Also, we study the same problems on K\"ahler foliations with a complete bundle-like metric.

\section{Introduction}
Let $(M,g_M,\mathcal F)$ be a foliated Riemannian manifold  with a foliation $\mathcal F$ of codimension $q$ and a bundle-like metric $g_M$ with respect to $\mathcal F$. 
  A
transversal conformal Killing field is a normal field with a flow
preserving the conformal class of the transverse metric. As a generalization of a transversal conformal Killing field, we define  the {\it transverse conformal Killing $r$-forms} $\phi$ as follows: for any
vector field $X$ normal to the foliation,
\begin{align*}
\nabla_X\phi={1\over r+1}i(X)d\phi-{1\over
q-r+1}X^{\frak b}\wedge\delta_T\phi,
\end{align*}
where  $X^{\frak b}$ is the dual 1-form
of $X$. For the definition of $\delta_T$,  see Section 2. The transverse conformal Killing
form $\phi$ with $\delta_T\phi=0$ is called {\it transverse Killing form}, which is a generalization of a transversal Killing field. There are many researches about such fields [\ref{AY},\ref{KT2},\ref{NT}]. Moreover, transverse (conformal) Killing forms on a point foliation are just (conformal) Killing forms, which were studied by many authors [\ref{KT},\ref{KA},\ref{MS},\ref{SE},\ref{TA}].
Recently, the  transverse Killing and conformal Killing forms were studied in [\ref{JU1},\ref{JU2},\ref{JJ2},\ref{JK}]. In particular, the non-existence of the transverse Killing and conformal Killing forms on compact foliated Riemannian manifolds was studied by S. D. Jung and K. Richardson in 2012 [\ref{JK}].  And the properties of such forms on K\"ahler foliations were studied by S. D. Jung and M. J. Jung  (for the transverse Killing forms) in 2012 [\ref{JJ2}] and  by S. D. Jung (for the transverse conformal Killing forms) in 2015 [\ref{JU2}], respectively. 

In 2014, S. D. Jung [\ref{JU1}]  studied the non-existence of the $L^2$-transverse Killing forms on complete foliated Riemannian manifolds. 

In this paper, we study the $L^2$-transverse conformal Killing forms for Riemannian  and  K\"ahler foliations with a bundle-like metric. In fact, we prove the following theorem.

\bigskip
\noindent {\bf Theorem A.} (cf. Theorem 3.6) {\it Let $(M,g_M,\mathcal F)$ be a complete foliated Riemannian manifold whose all leaves are compact. Assume that the mean curvature form $\kappa_B$ is bounded and coclosed. If the curvature endomorphism $F$ is nonpositive, then every $L^2$-transverse conformal Killing $r (1\leq r \leq q-1)$-form is parallel.}

\bigskip
\noindent{\bf Theorem B.} (cf. Corollary 3.7) {\it Let $(M,g_M,\mathcal F)$ be as in Theorem A. Assume that $\kappa_B$ is bounded and coclosed. If the curvature endomorphism $F$ is nonpositive and either negative at some point or ${\rm Vol}(M)=\infty$, then every $L^2$-transverse conformal Killing $r (1\leq r\leq q-1)$-form is trivial.}

\bigskip

\noindent Since the transversal conformal Killing field is a dual vector field of the transverse conformal Killing $1$-form, we have the following theorem.

\bigskip
\noindent {\bf Theorem C.} (cf. Corollary 3.9) {\it Let $(M,g_M,\mathcal F)$ be as in Theorem A. Assume that $\kappa_B$ is bounded and coclosed. If the transversal Ricci curvature is nonpositive and either negative at some point or ${\rm Vol}(M)=\infty$, then every $L^2$-transversal conformal Killing field is trivial.}

\bigskip
\noindent {\bf Remark.} Theorem C was proved by S. Yorozu [\ref{YO2}] when $\mathcal F$ is point foliation, and  by T. Aoki and S. Yorozu [\ref{AY}], by S. Nishikawa and Ph. Tondeur [\ref{NT1}]  when $\mathcal F$ is a minimal foliation. 

\bigskip
\noindent
Let $J$ be the extension of complex structure $J$ to basic forms on K\"ahler foliations. For details, see (4.3) in section 4. Then we have the following theorem.

\bigskip
\noindent{\bf Theorem D.} (cf. Theorem 4.12) {\it Let $(M,g_M,\mathcal F,J)$ be a complete Riemannian manifold with a K\"ahler foliation of codimension $q=2m>4$ whose all leaves are compact. If $\mathcal F$ is minimal, then for a $L^2$-transverse conformal Killing $r (2\leq r \leq q-2)$-form $\phi$, $J\phi$ is parallel.}

 \section{Transverse conformal Killing forms}
 Let $(M,g_M,\mathcal F)$ be a $(p+q)$-dimensional
Riemannian manifold with a foliation $\mathcal F$ of codimension
$q$ and a bundle-like metric $g_M$ with respect to $\mathcal F$.
 Then there exists an exact sequence of vector bundles
\begin{align}\label{eq1-1}
 0 \longrightarrow T\mathcal F \longrightarrow
TM {\overset\pi\longrightarrow} Q \longrightarrow 0,
\end{align}
where $T\mathcal F$ is the tangent bundle  and $Q=TM/T\mathcal F$ is the normal bundle
of $\mathcal F$. The metric $g_M$ determines an orthogonal
decomposition $TM=T\mathcal F\oplus T\mathcal F^\perp$, identifying $Q$ with $T\mathcal F^\perp$
and inducing a metric $g_Q$ on $Q$. 
 Let $\Omega_B^*(\mathcal F)$ be
the space of all {\it basic forms} on $M$, i.e.,
\begin{align}\label{eq1-4}
\Omega_B^*(\mathcal F)=\{\phi\in\Omega^*(M)\ | \ i(X)\phi=0,\
i(X)d\phi=0,
                     \quad \forall  X\in T\mathcal F\}.
\end{align}
Then  
$\Omega^*(M)=\Omega_B^*(\mathcal F) \oplus
\Omega_B^*(\mathcal F)^\perp [\ref{LO}].
$
Let $\nabla $ be the transverse Levi-Civita connection on $Q$ [\ref{KT}], which is extended to $\Omega _{B}^*(\mathcal F)$.
  The exterior differential $d$ on the  de Rham complex
$\Omega^*(M)$ restricts a differential $d_B:\Omega_B^r(\mathcal
F)\to \Omega_B^{r+1}(\mathcal F)$.  
    Let $\kappa\in Q^*$ be the
mean curvature form of $\mathcal F$. It is well-known
 that the basic part $\kappa_B$ of $\kappa$ is closed [\ref{LO}].
 We now recall the star operator $\bar
*:\Omega^r_B(\mathcal F)\to \Omega_B^{q-r}(\mathcal F)$, which is defined by
[\ref{PR},\ref{TO}]
\begin{align}\label{eq1-6}
 \bar * \phi =(-1)^{p(q-r)}*(\phi\wedge \chi_{\mathcal
F}),\quad\forall \phi\in\Omega_B^r(\mathcal F),
\end{align}
where $\chi_{\mathcal F}$ is the characteristic form of $\mathcal
F$ and $*$ is the Hodge star operator associated to $g_M$. Then, for any $\phi,\psi\in\Omega_B^r(\mathcal F)$, it is well-known that  $\phi\wedge\bar
*\psi =\psi\wedge\bar *\phi$ and  ${\bar
*}^2\phi=(-1)^{r(q-r)}\phi$  [\ref{PR}]. Let $\nu$ be the transversal volume
form, i.e., $*\nu=\chi_\mathcal F$.  The pointwise inner product
$\langle\ , \ \rangle$ on $\Lambda^r Q^*$ is defined uniquely by
\begin{align}\label{eq1-8}
\langle\phi,\psi\rangle\nu=\phi\wedge \bar *\psi.
\end{align}
The global inner product $\ll\cdot,\cdot\gg_B$ on $L^2\Omega_B^r(\mathcal F)$ is given by
\begin{align}
\ll\phi,\psi\gg_B =\int_M \langle\phi,\psi\rangle\mu_M,
\end{align}
where $\mu_M=\nu\wedge\chi_{\mathcal F}$ is the volume form with respect to $g_M$. With respect to this scalar product, the formal adjoint operator  $\delta_B:\Omega_B^r(\mathcal F)\to \Omega_B^{r-1}(\mathcal F)$ of $d_B$ is given by 
 \begin{align}\label{eq1-10}
 \delta_B\phi=(-1)^{q(r+1)+1}\bar *d_T\bar *
 \phi=\delta_T\phi+ (-1)^{q(r+1)}\bar *(\kappa_B \wedge\bar *\phi),
 \end{align}
 where $d_T=d-\kappa_B\wedge$ and $\delta_T=(-1)^{q(r+1)+1}\bar
 *d\bar *$ is the formal adjoint operator of $d_T$ with respect to this scalar product.
  Trivially, for any basic form $\phi\in L^2\Omega_{B}^r(\mathcal F)$,
 \begin{align}
 i(\kappa_B^\sharp)\phi = (-1)^{q(r+1)}\bar *(\kappa_B \wedge\bar *\phi),
 \end{align}
 where $(\cdot)^\sharp$ is a $g_Q$-dual vector to $(\cdot)$ [\ref{JK}]. 
 The basic
Laplacian $\Delta_B$ is given by
 $ \Delta_B = d_B\delta_B
+ \delta_B d_B$.
 Let $\{E_a\} (a=1,\cdots,q)$ be a local orthonormal basic frame on $Q$.
We define $ \nabla_{\rm tr}^*\nabla_{\rm tr} :\Omega_B^r(\mathcal F)\to
\Omega_B^r(\mathcal F)$ by
\begin{align}\label{eq1-12}
\nabla_{\rm tr}^*\nabla_{\rm tr}\phi =-\sum_a \nabla^2_{E_a,E_a}
\phi+\nabla_{\kappa_B^\sharp}\phi,\quad\phi\in\Omega_B^r(\mathcal F),
\end{align}
where $\nabla^2_{X,Y}=\nabla_X\nabla_Y -\nabla_{\nabla^M_XY}$ for
any $X,Y\in TM$ and $\nabla^M$ is the Levi-Civita connection on $M$. Then the operator $\nabla_{\rm tr}^*\nabla_{\rm tr}$
is positive definite and formally self adjoint on  $L^2\Omega_B^r(\mathcal F)$ [\ref{JU}]. We define the bundle map $A_Y:\Lambda^r
Q^*\to\Lambda^r Q^*$ for any $Y\in TM$ [\ref{KT2}]
by
\begin{align}\label{eq1-13}
A_Y\phi =\theta(Y)\phi-\nabla_Y\phi,
\end{align}
where $\theta(Y)$ is the transverse Lie derivative.
For any vector field $X\in T\mathcal F$, $\theta(X)\phi=\nabla_X\phi$  [\ref{KT1}] and so $A_X\phi=0$.
Now we define the curvature endomorphism $F:\Omega_B^r (\mathcal F)\to \Omega_B^r (\mathcal F)$ by 
\begin{align}
F(\phi)=\sum_{a,b}\theta^a \wedge i(E_b)R^\nabla(E_b,
 E_a)\phi,
 \end{align}
  where $R^\nabla$ is the curvature tensor with respect to $\nabla$ and $\theta^a$ is a dual 1-form to $E_a$. Note that if $\phi$ is a basic 1-form, then $F(\phi)^\sharp
 ={\rm Ric}^Q(\phi^\sharp)$, where ${\rm Ric}^Q$ is the transversal Ricci curvature of $\mathcal F$.
Now we recall 
the generalized Weitzenb\"ock formula.
\begin{thm} $[\ref{JU}]$ On a Riemannian foliation $\mathcal F$, we have that for any $\phi\in\Omega_B^r(\mathcal
  F)$,
\begin{align*}
  \Delta_B \phi = \nabla_{\rm tr}^*\nabla_{\rm tr}\phi +
  F(\phi)+A_{\kappa_B^\sharp}\phi.
\end{align*}
\end{thm}
From Theorem 2.1, we obtain that for any $\phi\in\Omega_B^r(\mathcal F)$,
\begin{align}\label{eq2-12}
\frac12\Delta_B|\phi|^2 = \langle\Delta_B\phi,\phi\rangle -|\nabla_{\rm tr}\phi|^2 -\langle F(\phi),\phi\rangle -\langle A_{\kappa_B^\sharp}\phi,\phi\rangle.
\end{align}
\begin{defn} {\rm A basic $r$-form $\phi\in \Omega_B^r(\mathcal F)$ is called a}  transverse conformal
Killing $r$-form {\rm if for any normal vector field $X\in T\mathcal F^\perp$,
\begin{align*}
\nabla_X\phi ={1\over r+1}i(X)d_B\phi -{1\over r^*+1}X^{\frak b}\wedge
\delta_T\phi,
\end{align*}
where $r^* = q-r$ and  $X^{\frak b}$ is the dual 1-form of $X$. In addition, if the
basic $r$-form $\phi$ satisfies $\delta_T\phi=0$, it is called a}
transverse Killing $r$-form.
\end{defn}
Note that a transverse conformal Killing 1-form (resp. transverse
Killing 1-form) is a $g_Q$-dual form of a transversal conformal
Killing  field (resp. transversal Killing  field) [\ref{JK}].

 \begin{prop}$[\ref{JK}]$ Let $\phi$ be a transverse conformal Killing $r$-form. Then
 \begin{align}
 & F(\phi) = {r\over r+1}\delta_Td_B\phi+{r^*\over r^*+1}d_B\delta_T\phi,\\
 &\nabla_{\rm tr}^*\nabla_{\rm tr}\phi={1\over r+1}\delta_Bd_B\phi +{1\over r^*+1}d_T\delta_T\phi.
 \end{align}
 \end{prop}
 \begin{thm} $[\ref{JK}]$ Any basic $r$-form $\phi$ is a transverse conformal Killing $r$-form if and only if $\bar *\phi$ is a transverse conformal Killing $(q-r)$-form.
 \end{thm}
   Now we recall the generalized maximum principle  on a complete foliated Riemannian manifold, that is, Riemannian foliation with a complete bundle-like metric.
\begin{thm} $[\ref{JJ}]$ 
Let $(M,g_M,\mathcal F)$ be a complete foliated Riemannian manifold whose all leaves are compact. Assume that $\kappa_B$ is bounded and coclosed.   Then a nonnegative basic function $f$ such that $(\Delta_B-\kappa_B^\sharp)f\leq 0$ with $\int_M f^p <\infty\ (p>1)$ is constant.
\end{thm}

 \section{ Vanishing theorems on Riemannian foliations  }
 
 Now we recall the vanishing theorem on a compact foliated Riemannian manifold.  
 \begin{thm}\label{prop2-11} $[\ref{JK}]$ Let $(M,g_M,\mathcal F)$ be a compact foliated Riemannian manifold with a foliation $\mathcal F$ of codimension $q$ and a bundle-like metric $g_M$
 such that $\delta_B\kappa_B=0$.
  Suppose $F$ is non-positive and negative at some point. Then, for any $1\leq r \leq q-1$, there are no non-trivial transverse
conformal Killing $r$-forms on $M$.
\end{thm}
In this section, we study some vanishing theorems of the $L^2$-transverse conformal Killing forms on complete foliated Riemannian manifold. The basic form $\phi$ is said to be $L^2$-{\it basic form} if $\phi\in L^2\Omega_B^* (\mathcal F)$, i.e.,  $\Vert\phi\Vert_B^2<\infty$.

Let $(M,g_M,\mathcal F)$ be a complete foliated Riemanian manifold and all leaves be compact. We consider a smooth function $\mu$ on $\mathbb R$ satisfying
\begin{align*}
(i) \ 0\leq \mu(t)\leq 1\ {\rm on}\ \mathbb R,\quad
(ii)\ \mu(t)=1\ \ {\rm for}\ t\leq 1,\quad
(iii)\ \mu(t)=0\ \ {\rm for}\ t\geq 2.
\end{align*}
Now, we fix a point $x_0\in M$. For each point $y \in M$, we
denote by $\rho(y)$ the distance between leaves through $x_0$ and $y$. For any real number $l>0$, we define a Lipschitz continuous function $\omega_l$  on $M$ by  
\begin{align*}
\omega_l(y)=\mu(\rho(y)/l).
\end{align*} 
Trivially, $\omega_l$ is a basic function. 
   Let $B(l) = \{ y \in M | \rho(y) \leq l \}$. Then 
$0 \leq \omega_l(y) \leq 1\  \text{for any} \ y \in M,\
\text{supp}\ \omega_l \subset B(2l),\  
\omega_l(y)=1 \ \text{for any} \  y \in B(l),\
\lim_{l \rightarrow \infty} \omega_l = 1 \ \text{and}\
|d_B\omega_l | \leq \frac{C}{l} \ \text{almost everywhere on}\ M,$
where $C$ is a positive constant independent of $l$ [\ref{YO1}]. Hence $\omega_l\psi$ has compact support for any basic form $\psi\in\Omega_B^*(\mathcal F)$ and $\omega_l\psi \to \psi$ (strongly) when $l\to \infty$.
\begin{lemma} $[\ref{KI}]$ For any $\phi\in\Omega_B^r (\mathcal F)$,  there exists a number $A$ depending only on $\mu$, such that
\begin{align*}
&\Vert d_B\omega_l \wedge\phi \Vert_{B(2l)}^2\leq {qA^2 \over l^2} \Vert\phi \Vert_{B(2l)}^2,\\
&\Vert d_B\omega_l \otimes \phi\Vert_{B(2l)}^2 \leq {qA^2\over l^2}\Vert\phi\Vert_{B(2l)}^2,
\end{align*}
where $\Vert \phi\Vert_{B(2l)}^2 =\int_{B(2l)}\langle\phi,\phi\rangle\mu_M$.
\end{lemma}
\begin{lemma} 
Let $(M,g_M,\mathcal F)$ be a complete foliated Riemannian manifold whose all leaves are compact. Assume that $\kappa_B$ is bounded and coclosed.  Then a $L^2$- transverse conformal Killing $r$-form $\phi$ satisfies
\begin{align*}
{1\over r+1} \lim_{l\to\infty} \ll i(\kappa_B^\sharp)d_B\phi,\omega_l^2\phi\gg_{B(2l)} = {1\over r^* +1}\lim_{l\to\infty}\ll\kappa_B\wedge\delta_T\phi,\omega_l^2\phi\gg_{B(2l)}.
\end{align*} 
\end{lemma}
{\bf Proof.} Let $\phi$ be a $L^2$-transverse conformal Killing $r$-form. By a direct calculation, we have
\begin{align*}
\kappa_B^{\sharp}(|\omega_l\phi|^2 ) = 2\langle \nabla_{\kappa_B^\sharp}(\omega_l\phi),\omega_l\phi\rangle = 2\langle \omega_l\kappa_B,d_B\omega_l\rangle |\phi|^2 +2\langle\nabla_{\kappa_B^\sharp}\phi,\omega_l^2\phi\rangle.
\end{align*} 
Since $\delta_B\kappa_B =0$, $\int_{B(2l)}\kappa_B^{\sharp}(|\omega_l\phi|^2)=0$. Hence we have
\begin{align*}
0=\int_{B(2l)} \langle\omega_l\kappa_B, d_B\omega_l\rangle |\phi|^2 + \int_{B(2l)} \langle\nabla_{\kappa_B^\sharp}\phi,\omega_l^2\phi\rangle.
 \end{align*}
 Since $|d\omega_l|<{C\over l}$ and $|\kappa_B|<\infty$, 
 $\lim_{l\to\infty} \int_M \langle\omega_l\kappa_B,d_B\omega_l\rangle|\phi|^2 = 0$. Hence  
 \begin{align}
 \lim_{l\to\infty} \int_{B(2l)}\langle\nabla_{\kappa_B^\sharp} \phi,\omega_l^2\phi\rangle =0.
 \end{align} 
 From Definition 2.2, the proof follows. $\Box$

\begin{lemma} Let $(M,g_M,\mathcal F)$ be as in Lemma 3.3.  If $\phi$ is a $L^2$-basic form, then 
\begin{align}
&\lim_{l\to\infty}\ll \omega_l d_B\phi,d_B\omega_l\wedge\phi\gg_{B(2l)}\ \geq - A_1 \Vert d_B\phi\Vert_B^2,\\
&\lim_{l\to\infty}\ll\omega_l\delta_T\phi,i(\nabla \omega_l)\phi\gg_{B(2l)}\ \geq -A_2 \Vert\delta_T\phi\Vert_B^2,\\
&\lim_{l\to\infty}\ll\kappa_B\wedge\delta_T\phi,\omega_l^2\phi\gg_{B(2l)}\ \geq -A_3  \Vert\delta_T\phi\Vert_B^2 -{1\over 4A_3} \Vert i(\kappa_B^\sharp)\phi\Vert_B^2
\end{align}
for any positive real numbers $A_1, A_2$ and $A_3$.
\end{lemma}
{\bf Proof.} Let $\phi$ be a $L^2$-basic form. From Lemma 3.2 and the Schwartz inequality, we have
\begin{align*}
2\Big|\ll \omega_l d_B\phi,d_B\omega_l\wedge\phi\gg_{B(2l)}\Big| \leq {\epsilon_1}\Vert\omega_l d_B\phi\Vert_{B(2l)}^2 + {qA^2\over l^2\epsilon_1}\Vert\phi\Vert_{B(2l)}^2
\end{align*}
for a positive real number $\epsilon_1$.
If we let $l\to \infty$, then (3.2) is proved. Similarly, by using the inequality $|d_B\omega_l|^2|\phi|^2 = |i(\nabla\omega_l)\phi|^2 + |d_B\omega_l\wedge\phi|^2$, we have
\begin{align*}
2\Big|\ll \omega_l \delta_T\phi,i(\nabla\omega_l)\phi\gg_{B(2l)}\Big| &\leq {\epsilon_2}\Vert\omega_l \delta_T\phi\Vert_{B(2l)}^2 + {C^2-qA^2\over l^2\epsilon_2}\Vert\phi\Vert_{B(2l)}^2
\end{align*}
for a positive real number $\epsilon_2$. If we let $l\to\infty$, then (3.3) is proved. From the Schwartz's inequality, the proof of (3.4) is trivial.  $\Box$

From Lemma 3.3 and Lemma 3.4, we have the following proposition.
\begin{prop} Let $(M,g_M,\mathcal F)$ be as in Lemma 3.3. Assume that  $\kappa_B$ is bounded and coclosed. Then
\begin{align*}
&\limsup_{l\to\infty}\ll F(\phi),\omega_l^2\phi\gg_{B(2l)}\\
&\geq {r\over {r+1}}(1-2A_1)\Vert d_B\phi\Vert_B^2 + {1\over {r^*+1}}\{r^*(1-2A_2)-A_3(q-2r)\}\Vert\delta_T\phi\Vert_B^2\\
& -{q-2r\over 4(r^*+1)A_3} \Vert i(\kappa_B^\sharp)\phi\Vert_B^2 
\end{align*}
for any positive real numbers $A_1,A_2$ and $A_3$.
\end{prop}
{\bf Proof.} Let $\phi$ be a $L^2$-transverse conformal Killing $r$-form. Since $\delta_T(\omega_l^2\phi)=\omega_l^2\delta_T\phi -2\omega_l i(\nabla\omega_l)\phi$,  from (2.12) we have 
\begin{align*}
&\ll F(\phi),\omega_l^2\phi\gg_{B(2l)}\\
 &={r\over r+1}\Vert\omega_l d_B\phi\Vert_{B(2l)}^2 +{r^*\over r^*+1}\Vert\omega_l \delta_T\phi\Vert_{B(2l)}^2\\
&+{2r\over r+1}\ll\omega_l d_B\phi,d_B\omega_l\wedge\phi\gg_{B(2l)} -{2r^*\over r^*+1}\ll\omega_l\delta_T\phi,i(\nabla\omega_l)\phi\gg_{B(2l)}\\
&-{r\over r+1}\ll i(\kappa_B^\sharp)d_B\phi,\omega_l^2\phi\gg_{B(2l)} +{r^*\over r^*+1}\ll\kappa_B\wedge \delta_T\phi,\omega_l^2\phi\gg_{B(2l)}.
\end{align*}
From Lemma 3.3, we have
\begin{align*}
\limsup_{l\to\infty}\ll F(\phi),\omega_l^2\phi\gg_{B(2l)}
 &={r\over r+1}\Vert d_B\phi\Vert_{B(2l)}^2 +{r^*\over r^*+1}\Vert \delta_T\phi\Vert_{B(2l)}^2\\
&+{2r\over r+1}\limsup_{l\to\infty}\ll\omega_l d_B\phi,d_B\omega_l\wedge\phi\gg_{B(2l)} \\
&-{2r^*\over r^*+1}\limsup_{l\to\infty}\ll\omega_l\delta_T\phi,i(\nabla\omega_l)\phi\gg_{B(2l)}\\
&+{q-2r\over r^*+1}\limsup_{l\to\infty}\ll\kappa_B\wedge \delta_T\phi,\omega_l^2\phi\gg_{B(2l)}.
\end{align*}
From Lemma 3.4, the proof is completed. $\Box$

Hence we have the following theorem.
\begin{thm} Let $(M,g_M,\mathcal F)$ be as in Lemma 3.3.  Assume that $\kappa_B$ is bounded and coclosed.  If the curvature endomorphism $F$ is nonpositive, then every $L^2$-transverse conformal Killing $r  (1\leq r\leq q-1)$-form is  parallel. 
\end{thm}
{\bf Proof.} Let $\phi$ be a $L^2$-transverse conformal Killing $r$-form. 

(i) In case of $r\geq {q\over 2}$. From Proposition 3.5, we have
\begin{align*}
\limsup_{l\to\infty}\ll F(\phi),\omega_l^2\phi\gg_{B(2l)}
&\geq {r\over {r+1}}(1-2A_1)\Vert d_B\phi\Vert_B^2\\
& + {1\over {r^*+1}}\{r^*(1-2A_2)-A_3(q-2r)\}\Vert\delta_T\phi\Vert_B^2.
\end{align*}
If we choose $0<A_1<\frac12$ and $0<A_2<\frac12$, then
\begin{align}
\limsup_{l\to\infty}\ll F(\phi),\omega_l^2\phi\gg_{B(2l)}\geq 0.
\end{align}
Since $F$ is nonpositive, from (3.5) we have
\begin{align}
\limsup_{l\to\infty}\ll F(\phi),\omega_l^2\phi\gg_{B(2l)}= 0. 
\end{align}
Hence from Proposition 3.5, we have $d_B\phi=0, \ \delta_T\phi=0$ and $i(\kappa_B^\sharp)\phi=0$.  Also from (2.13), we have
\begin{align}
\nabla_{\rm tr}^*\nabla_{\rm tr}\phi=0.
\end{align}
By multiplying $\omega_l^2\phi$ and by the Schwarz's inequality, we have that, for a real number $\epsilon>0$
\begin{align*}
0&=\ll\nabla_{\rm tr}^*\nabla_{\rm tr}\phi,\omega_l^2\phi\gg_{B(2l)} \\
&=\Vert\omega_l\nabla_{\rm tr}\phi\Vert_{B(2l)}^2 +\ll \nabla_{\rm tr}\phi, 2\omega_l d_B\omega_l\otimes \phi\gg_{B(2l)} \\
&\geq (1-\epsilon)\Vert\omega_l\nabla_{\rm tr}\phi\Vert_{B(2l)}^2 -{1\over\epsilon}\Vert d_B\omega_l\otimes\phi\Vert_{B(2l)}^2. 
\end{align*}
From Lemma 3.2, we have
\begin{align*}
0\geq (1-\epsilon)\Vert\omega_l\nabla_{\rm tr}\phi\Vert_{B(2l)}^2 -{qA^2\over l^2\epsilon}\Vert\phi\Vert_{B(2l)}^2.
\end{align*}
If we let $l\to\infty$, then $(1-\epsilon)\Vert \nabla_{\rm tr}\phi \Vert_B^2\leq 0$. So if we choose $\epsilon<1$, then
\begin{align}
\Vert \nabla_{\rm tr}\phi \Vert_B^2 =0.
\end{align}
Hence $\nabla_{\rm tr}\phi=0$, that is, $\phi$ is parallel. 

(ii) In case of $r\leq {q\over 2}$.  From Theorem 2.4, $\bar *\phi$ is also $L^2$-transverse conformal Killing $r^*$-form. Since $r^*\geq {q\over 2}$, from (i), $\bar *\phi$ is parallel. Since $\nabla_{\rm tr}\bar*\phi =\bar *\nabla_{\rm tr}\phi$ and $\bar *$ is an isometry, $\phi$ is parallel. Hence any  $L^2$-transverse conformal Killing $r$-form is parallel. From (i) and (ii), the proof is completed. $\Box$
\begin{coro} Let $(M,g_M,\mathcal F)$ be as in Lemma 3.3.  Assume that $\kappa_B$ is bounded and coclosed.
If $F$ is nonpositive  and either negative at some point or ${\rm Vol}(M)=\infty$, then every $L^2$-transverse conformal Killing $r\ (1\leq r\leq q-1)$-form $\phi$  is trivial.
\end{coro}
{\bf Proof.} Since $\phi$ is parallel, $F(\phi)=0$. Hence the negativity of $F$ means that $\phi$ is trivial. Now, we consider ${\rm Vol}(M)=\infty$. From (2.11) and Theorem 2.5, $|\phi|^2$ is constant. Hence $\int_M |\phi|^2<\infty$ and ${\rm Vol}(M)=\infty$ yield that $\phi$ is trivial. $\Box$
\begin{rem}  {\rm $(cf.\ [\ref{JU1}])$ Let $(M,g_M,\mathcal F)$ be as in Lemma 3.3. Assume that $\kappa_B$ is bounded and coclosed. If $F$ is nonpositive and either negative at some point or ${\rm Vol}(M)=\infty$, then every $L^2$-transverse Killing $r\ (1\leq r\leq q-1)$-form is trivial.}
\end{rem}
Since $F(\phi)^\sharp = {\rm Ric}^Q(\phi^\sharp)$ for any basic $1$-form, we have the following corollary.
\begin{coro} Let $(M,g_M,\mathcal F)$ be as in Lemma 3.3. Assume that $\kappa_B$ is bounded and coclosed. If the transversal Ricci curvature is nonpositive and either negative at some point or ${\rm Vol}(M)=\infty$, then every $L^2$-transversal conformal Killing field is trivial.
\end{coro}

\begin{rem} {\rm  (1) When $\mathcal F$ is a foliation by points, Corollary 3.9 was given in [\ref{YO2}].

(2) When $\mathcal F$ is minimal, Corollary 3.9 was proved in [\ref{AY}] and [\ref{NT1}], respectively. So Corollary 3.9 is a generalization of the results in [\ref{AY}, \ref{NT1}] to the non-minimal case. 
}
\end{rem}

\section{ The properties on  K\"ahler foliations}
Let $(M,g_M,\mathcal F,J)$ be a Riemannian manifold with a
K\"ahler foliation $\mathcal F$ of codimension $q=2m$ and a
bundle-like metric $g_M$ [\ref{JU2},\ref{NT}]. Namely, there is a holonomy invariant almost complex structure $J:Q\to Q$ with respect to which $g_Q$ is Hermitian, i.e., $g_Q(JX,JY)=g_Q(X,Y)$ for any $X,Y\in Q$ and $\nabla J=0$. Note that for any $X,
Y\in TM$,
\begin{equation}
\Omega(X,Y)=g_Q(\pi(X),J\pi(Y))
\end{equation}
defines a basic 2-form $\Omega$, which is closed as consequence of
$\nabla g_Q=0$ and $\nabla J=0$. 
 Now, we define the operators $L:\Omega_B^r(\mathcal F)\to \Omega_B^{r+2}(\mathcal F)$ and $\Lambda:\Omega_b^r(\mathcal F)\to \Omega_B^{r-2}(\mathcal F)$ respectively by $[\ref{JJ3}]$
\begin{align}
L(\phi)=\epsilon(\Omega)\phi,\quad\Lambda(\phi)=i(\Omega)\phi,
\end{align}
where $\epsilon(\Omega)\phi=\Omega\wedge\phi$ and $i(\Omega)=-\frac12 \sum_{a=1}^{2m} i(JE_a)i(E_a)$. Trivially, for any basic forms $\phi\in\Omega_B^r(\mathcal F)$ and $\psi\in\Omega_B^{r+2}(\mathcal F)$, $\langle L(\phi),\psi\rangle =\langle\phi,\Lambda(\psi)\rangle$. Moreover, for any basic $r$-form $\phi$, $[\Lambda,L]\phi = \frac12(q-2r)\phi$.

Also, we define the operator
$\tilde J:\Omega_B^r(\mathcal F)\to \Omega_B^r(\mathcal F)$ by
\begin{align} 
\tilde J(\phi) =\sum_{a=1}^{2m}J\theta^a\wedge i(E_a)\phi.
\end{align}
Trivially,  $\tilde J\phi =J\phi$ for any basic 1-form $\phi$. So $\tilde J$ is an extension of  the complex structure $J$ to basic forms.  From now on, if we have no confusion, we write $\tilde J\equiv J$. 

  \begin{lemma} $[\ref{JU2}]$ On a K\"ahler foliation $(\mathcal F,J)$,  we have
\begin{gather*}
[J,L]=[J,\Lambda]=[F,J]=[F,\Lambda]=0.
\end{gather*}
In particular, if $\mathcal F$ is minimal, then
\begin{gather*}
 [J,\Delta_B]=[\Lambda,\Delta_B]=0.
 \end{gather*}
\end{lemma}
 \begin{prop} $[\ref{JU2}]$ On a K\"ahler foliation $(\mathcal F,J)$, a transverse conformal Killing $r$-form $\phi$ satisfies
 \begin{align*}
 F(J\phi)=0.
 \end{align*}
 \end{prop}
Now, we recall the operators $d_B^c:\Omega_B^r(\mathcal F)\to\Omega_B^{r+1}(\mathcal F)$ and $\delta_B^c:\Omega_B^r(\mathcal F)\to\Omega_B^{r-1}(\mathcal F)$, which are given by $[\ref{JJ3}]$
\begin{align}
&d_B^c\phi=\sum_{a=1}^{2m}J\theta^a\wedge\nabla_{E_a}\phi,\\ 
&\delta_B^c \phi= -\sum_{a=1}^{2m}i(JE_a)\nabla_{E_a}\phi + i(J\kappa_B^\sharp)\phi.
\end{align}
 Trivially, $\delta_B^c$ is a formal adjoint of $d_B^c$ and ${\delta_B^c}^2 = {d_B^c}^2 =0$ [\ref{JJ3}]. Also, we define two operators $d_T^c$ and $\delta_T^c$ by
 \begin{align}
 d_T^c = d_B^c -\epsilon(J\kappa_B),\quad \delta_T^c = \delta_B^c -i(J\kappa_B^\sharp).
 \end{align}
 If $\mathcal F$ is minimal, then $d_B^c =d_T^c$ and $\delta_B^c =\delta_T^c$.
Then we have the following lemma.
\begin{lemma} $[\ref{JU2}]$ On  a minimal K\"ahler foliation $(\mathcal F,J)$, we have
\begin{align}
d_B^c\delta_B + \delta_B d_B^c =d_B \delta_B^c +\delta_B^c d_B= 0.
\end{align}
\end{lemma}

 \begin{lemma} $[\ref{JU2}]$ On a K\"ahler foliation $(\mathcal F,J)$, a transverse conformal Killing $r$-form $\phi$ satisfies
 \begin{align}
 (rr^*-r-2)d_B^c\phi&=(r^*+1)d_BJ\phi -2(r+1)\delta_T L\phi,\\
 (rr^*-r^*-2)\delta_T^c\phi&=(r+1)\delta_TJ\phi +2(r^*+1)d_B\Lambda\phi.
 \end{align}
 \end{lemma}
 \begin{lemma} $[\ref{JU2}]$ Let $(M,g_M,J,\mathcal F)$ be a K\"ahler foliation. Then a transverse conformal Killing $r$-form $\phi$ satisfies
\begin{align*}
a_1 \delta_T d_B J\phi + a_2 d_B\delta_T J\phi + a_3e(\kappa_B^\sharp)L\phi=0,
\end{align*}
where $a_1 = (r^*+1)(rr^*-r^*-r)(rr^*-r^*-2),\ a_2 = (r^*+1)(rr^*-r^*-r)(rr^*-r-2)$ and $a_3={2(r-r^*)(rr^*-r^*-2)}$.
\end{lemma}
Note that $rr^*-r^*-2=0$ if and only if $q=4$. Hence we have the following.
\begin{lemma} $[\ref{JU2}]$ Let $(M,g_M,J,\mathcal F)$ be a minimal K\"ahler foliation of codimension $q(\ne 4)$. Then for  a transverse conformal Killing $r$ $(2\leq r\leq q-2)$-form $\phi$,
\begin{align}
&b_1b_3 \delta_B d_B J\Lambda\phi = (1-b_1b_3)\delta_Bd_B^c\Lambda\phi + b_3(1-b_1)\delta_B^c d_B\Lambda\phi,\\
&b_2 d_B \delta_B J\Lambda\phi + b_2(1-b_1)\delta_B d_B J\Lambda\phi =(b_1b_2-1)\delta_B^c d_B \Lambda\phi + b_2(b_1-1)\delta_B d_B^c\Lambda\phi,
\end{align}
where $b_1= {r^*(r+1)\over rr^* -r^*-2},\ b_2={r+1\over (r^*+1)(r-1)}$ and $b_3={r^*+1\over (r+1)(r^*-1)}$.
\end{lemma}
From now on, let $(M,g_M,\mathcal F,J)$ be a complete K\"ahler foliation, i.e., a
K\"ahler foliation with a complete bundle-like metric.

\begin{thm} Let $(M,g_M,\mathcal F,J)$ be a complete K\"ahler foliation of codimension $q=2m$ whose all leaves are compact. Assume that $\kappa_B$ is bounded and coclosed. Let $\phi\in \Omega_B^m(\mathcal F)$ be a $L^2$-transverse conformal Killing ${q\over 2}$-form. If $q\ne 4$, then $J\phi\in \Omega_B^m(\mathcal F)$ is parallel. In addition, if $q=4$ and $\mathcal F$ is minimal, then $J\phi\in \Omega_B^2(\mathcal F)$ is parallel.
\end{thm}
{\bf Proof.} The proof is similar process to the one in [\ref{JU2}, Theorem 5.6].  Let $\phi\in \Omega_B^m(\mathcal F )$ be a $L^2$-transverse conformal Killing $m$-form. By Lemma 4.5, $a_1=a_2=m(m+1)^3 (m-2)^2$ and $a_3=0$. 

(i) In cse of $m\ne 2$, i.e., $q\ne 4$,  we have $a_1=a_2\ne 0$.  That is,
\begin{align*}
d_B\delta_TJ\phi + \delta_T d_B J\phi=0.
\end{align*}
Equivalently,  we have
\begin{align}
\Delta_B J\phi =\theta(\kappa_B^\sharp)J\phi.
\end{align}
Hence, by the scalar Weitzenb\"ock formula (2.11) and Proposition 4.2,  we have 
\begin{align}
\frac12(\Delta_B-\kappa_B^\sharp) |J\phi|^2 = - |\nabla_{\rm tr}J\phi|^2\leq 0.
\end{align} 
From the generalized maximum principle (Theorem 2.5), $|J\phi|$ is constant. Again, from (4.13), we have
\begin{align}
\nabla_{\rm tr}J\phi =0,
\end{align}
which implies that $J\phi$ is parallel. 

(ii) In case of $m=2$, i.e., $q=4$,  from Lemma 4.4, we  have
\begin{align}
d_BJ\phi =2\delta_TL\phi,\quad \delta_TJ\phi =-2d_B\Lambda\phi.
\end{align} 
From (4.15), we have
\begin{align*}
\delta_B d_B J\phi =2\delta_B\delta_TL\phi,\quad d_B\delta_TJ\phi =0.
\end{align*}
Hence if $\mathcal F$ is minimal, then $\Delta_B J\phi=0$.
From Theorem 2.1 and Proposition 4.2, we have
\begin{align}
\nabla_{\rm tr}^*\nabla_{\rm tr}J\phi =0.
\end{align} 
By multiplying $\omega_l^2J\phi$ and by integrating, we have
\begin{align}
\Vert\omega_l \nabla_{\rm tr} J\phi\Vert^2_{B(2l)} + 2\ll\omega_l\nabla_{\rm tr}J\phi,d_B\omega_l\otimes J\phi\gg_{B(2l)}=0.
\end{align}
By the Schwarz inequality and Lemma 3.2, we have
\begin{align*}
\lim_{l\to\infty} \ll \omega_l \nabla_{\rm tr} J\phi,d_B\omega_l\otimes J\phi\gg_{B(2l)} =0.
\end{align*}
Hence from (4.17), we have $\Vert\nabla_{\rm tr}J\phi\Vert_B^2 =0$, i.e., $J\phi$ is parallel.  Consequently, from (i) and (ii), the proof is completed. $\Box$

\begin{coro} Let $(M,g_M,\mathcal F,J)$ be as in Theorem 4.7, and suppose that $\mathcal F$ is minimal. Then for any $L^2$-transverse conformal Killing ${q\over 2}$-form, $J\phi$ is parallel.
\end{coro}
\begin{lemma} Let $(M,g_M,\mathcal F,J)$ be as in Theorem 4.7, suppose that $\mathcal F$ is minimal. Then for any $L^2$-basic form $\phi$, 
\begin{align}
&\lim_{l\to\infty}\ll \omega_l^2 \delta_B d_B J\Lambda\phi, \delta_B^c d_B \Lambda\phi\gg_{B(2l)} =0,\\
& \lim_{l\to\infty}  \ll \omega_l^2 \delta_B d_B^c\Lambda\phi,\delta_B^c d_B\Lambda\phi\gg_{B(2l)} = 0,\\
&\lim_{l\to\infty}\ll \omega_l^2\delta_B\{d_BJ\Lambda\phi+d_B^c\Lambda \phi\},d_B\delta_BJ\Lambda\phi\gg_{B(2l)}=0.
\end{align}
\end{lemma}
{\bf Proof.} 
Note that for any $\phi\in L^2 \Omega_B^r (\mathcal F)$, we have
\begin{align}
\omega_l^2 \delta_B d_B J\Lambda\phi &= \delta_B \{d_B(\omega_l^2 J\Lambda\phi) - 2\omega_l d_B\omega_l \wedge J\Lambda\phi\} + 2\omega_l i(\nabla\omega_l) d_B J\Lambda\phi,\\
\omega_l^2\delta_B d_B^c\Lambda\phi &=  \delta_B \{d_B^c(\omega_l^2\Lambda\phi)-2\omega_lJ(d_B\omega_l)\wedge\Lambda\phi\}+2\omega_l i(\nabla\omega_l) d_B^c\Lambda\phi.
\end{align}
Since $\mathcal F$ is minimal, from Lemma 4.3 and (4.21), we have
\begin{align}
\ll \omega_l^2 \delta_B d_B J\Lambda\phi, \delta_B^c d_B \Lambda\phi\gg_{B(2l)}
= -2\ll\omega_l d_B J\Lambda\phi,  d_B\omega_l\wedge d_B \delta_B^c \Lambda\phi\gg_{B(2l)}.
\end{align}
From (4.23),  by using the Schwarz inequality and Lemma 3.2, we have
\begin{align*}
\Big|\ll \omega_l^2 \delta_B d_B J\Lambda\phi, d_B\delta_B^c  \Lambda\phi\gg_{B(2l)}\Big|
&\leq \epsilon \Vert \omega_l d_B J\Lambda\phi\Vert_{B(2l)}^2 + {1\over\epsilon}\Vert d_B\omega_l\wedge d_B\delta_B^c \Lambda\phi\Vert_{B(2l)}^2\\
&\leq\epsilon \Vert \omega_l d_B J\Lambda\phi\Vert_{B(2l)}^2 + {qA^2\over\epsilon l^2}\Vert d_B\delta_B^c \Lambda\phi\Vert_{B(2l)}^2
\end{align*}
for any positive real number $\epsilon$. If we let  $l\to\infty$, then
\begin{align*}
\lim_{l\to\infty}\Big|\ll \omega_l^2 \delta_B d_B J\Lambda\phi, \delta_B^c d_B \Lambda\phi\gg_{B(2l)}\Big|
&\leq\epsilon \Vert d_B J\Lambda\phi\Vert_B^2.
\end{align*}
Since $\epsilon$ is arbitrary, the proof of (4.18) is completed. By using (4.22), the proof of (4.19) is similarly completed. Similarly, the proof of (4.20) follows.  $\Box$  


\begin{thm} Let $(M,g_M,\mathcal F,J)$ be as in Theorem 4.7,  and suppose that $\mathcal F$ is minimal.  Then for any $L^2$-transverse conformal Killing $r$ $(2\leq r\leq q-2)$-form $\phi$, $J\Lambda\phi$ is basic-harmonic.
\end{thm}
{\bf Proof.}
(i) In case of $q\ne 4$. From (4.10), we have 
\begin{align}
b_1b_3\ll \omega_l^2 \delta_B d_B J\Lambda\phi,\delta_B^c d_B\Lambda\phi\gg_{B(2l)}&=(1-b_1b_3)\ll\omega_l^2 \delta_B d_B^c\Lambda\phi,\delta_B^c d_B\Lambda\phi\gg_{B(2l)}\notag \\
&+ b_3(1-b_1)\Vert\omega_l \delta_B^c d_B\Lambda\phi\Vert^2_{B(2l)}.
\end{align}
Since $b_1\ne 1$ and $b_3\ne 0$, from Lemma 4.9, if we let $l\to\infty$, then
\begin{align}
\Vert \delta_B^c d_B\Lambda\phi\Vert_B^2 =0, \quad i.e., \ \delta_B^c d_B\Lambda\phi=0.
\end{align}
Hence from (4.11), we have
\begin{align}
b_2 d_B\delta_B J\Lambda\phi = b_2(b_1-1)\delta_B\{d_BJ\Lambda\phi+d_B^c\Lambda \phi\}.
\end{align}
By multiplying $\omega_l^2 d_B\delta_BJ\Lambda\phi$ in (4.26) and by integrating, we have
\begin{align*}
b_2\Vert\omega_l d_B\delta_BJ\Lambda\phi\Vert^2 =b_2(b_1-1)\ll \omega_l^2\delta_B\{d_BJ\Lambda\phi+d_B^c\Lambda \phi\},d_B\delta_BJ\Lambda\phi\gg_{B(2l)}.
\end{align*}
If we let $l\to\infty$, then from (4.20),
\begin{align}
\Vert d_B\delta_BJ\Lambda\phi\Vert_B^2 =0,\quad i.e., \ d_B\delta_B J\Lambda\phi=0.
\end{align}
Again from Lemma 4.6, since $b_2(1-b_1)\ne 0$, we have
\begin{align}
b_1 b_3 \delta_B d_B J\Lambda\phi &= (1-b_1 b_2) \delta_B d_B^c\Lambda\phi\\
 \delta_Bd_BJ\Lambda\phi &=-\delta_B d_B^c\Lambda\phi.
 \end{align}
Since $b_1b_3\ne b_1 b_2-1$, from (4.28) and (4.29), we have
\begin{align}
\delta_B d_B^c\Lambda\phi=0\quad {and}\quad \delta_B d_B J\Lambda\phi=0.
\end{align}  
From (4.27) and (4.30), we have
\begin{align}
\Delta_B J\Lambda\phi=0.
\end{align}
That is, $J\Lambda\phi$ is basic-harmonic. (ii) In case of $q=4$, since $\mathcal F$ is minimal, from Theorem 4.7, $J\phi\in\Omega_B^2(\mathcal F)$ is parallel and so basic-harmonic, i.e., $\Delta_BJ\phi=0$. Hence from Lemma 4.1, $\Delta_BJ\Lambda\phi=0$, i.e., $J\Lambda\phi$ is a basic-harmonic 2-form. Hence from (i) and (ii), the proof is completed.   $\Box$ 
    
\begin{coro}  Let $(M,g_M,\mathcal F,J)$ be as in Theorem 4.7, and suppose that $\mathcal F$ is minimal.  Then for any $L^2$-transverse conformal Killing $r$ $(2\leq r\leq q-2)$-form $\phi$, $J\Lambda\phi$ is parallel.
\end{coro}
{\bf Proof.} From Lemma 4.1 and Proposition 4.2, we have $F(J\Lambda\phi)=\Lambda F(J\phi)=0$.  Hence by the generalized Weitzenb\"ock formula (Theorem 2.1) and Theorem 4.10, we have
\begin{align}
\nabla_{\rm tr}^*\nabla_{\rm tr} J\Lambda\phi=0.
\end{align}
By multiplying $\omega_l^2 J\Lambda\phi$ in (4.32) and by integrating, we have
\begin{align}
\Vert\omega_l\nabla_{\rm tr}J\Lambda\phi\Vert_{B(2l)}^2 +2\ll \omega_l\nabla_{\rm tr}J\Lambda\phi,d_B\omega_l\otimes J\Lambda\phi\gg_{B(2l)}=0.
\end{align} 
By the Schwarz inequality and Lemma 3.2, we have
\begin{align}
\lim_{l\to\infty}\ll \omega_l\nabla_{\rm tr}J\Lambda\phi, d_B\omega_l\otimes J\Lambda\phi\gg_{B(2l)}=0.
\end{align}
Hence from (4.33) and (4.34), if we let $l\to\infty$, then 
\begin{align}
\nabla_{\rm tr} J\Lambda\phi=0.
\end{align}
That is, $J\Lambda\phi$ is parallel. $\Box$

\begin{thm}  Let $(M,g_M,J,\mathcal F)$ be as in Theorem 4.7, and suppose that $\mathcal F$ is minimal. Then for a  $L^2$-transverse conformal Killing $r$  $(2\leq r\leq q-2)$-form $\phi$,  $J\phi$ is parallel.
\end{thm}
{\bf Proof.} The proof is similar process to the one in [\ref{JU2}, Theorem 5.12]. That is, let $\phi$ be a $L^2$-transverse conformal Killing $r$-form. Then $\bar *\phi$ is also a $L^2$-transverse conformal Killing $(q-r)$-form [\ref{JK}]. Hence by Corollary 4.11, $J\Lambda  \bar *\phi$ is parallel. Since $[\nabla_{\rm tr},\bar *]=0$, $[J,\bar *] =0$ and $L\bar *=\bar *\Lambda$, $\bar * J\Lambda \bar * \phi = \pm LJ\phi$ is parallel.
Note that $(m-r)J\phi = [\Lambda,L]J\phi$. Since $[L,\nabla_{\rm tr}]=[\Lambda,\nabla_{\rm tr}]=0$ and $J\Lambda\phi$ is parallel, $(m-r)\nabla_{\rm tr}J\phi=0$. So if $r\ne m$, then $J\phi$ is parallel. For $r=m$, from Corollary 4.8, $J\phi$ is parallel. So the proof is completed. $\Box$

\bigskip
\noindent
For the point foliation, we have the following corollary.
\begin{coro} Let $(M,g_M, J)$ be a complete K\"ahler manifold of dimension $n=2m$. Then for any $L^2$-conformal Killing $r ( 2\leq r \leq n-2)$-forms $\phi$, $J\phi$ is parallel.
\end{coro}
\bigskip
\noindent{\bf Acknowledgements.} The first author was supported by the National Research Foundation of Korea(NRF) grant funded by the Korea government(MSIP) (NRF-2015R1A2A2A01003491) and the second author was supported by NSFC (No. 11371080).



\vskip 1.0cm 

\noindent Seoung Dal Jung

\noindent Department of Mathematics, Jeju National University, Jeju 690-756, Korea (e-mail: sdjung@jejunu.ac.kr)

\vskip 1.0cm

\noindent Huili Liu

\noindent Department of Mathematics, Northeastern University, Shenyang 110004, P. R. China (e-mail: liuhl@mail.neu.edu.cn)

\end{document}